\input amstex
\input amsppt.sty
\magnification=\magstep1
\hsize=30truecc
\vsize=22.2truecm
\baselineskip=16truept
\NoBlackBoxes
\TagsOnRight \pageno=1 \nologo
\def\Z{\Bbb Z}
\def\N{\Bbb N}

\def\C{\Bbb C}
\def\l{\left}
\def\r{\right}
\def\bg{\bigg}
\def\({\bg(}
\def\[{\bg\lfloor}
\def\){\bg)}
\def\]{\bg\rfloor}
\def\t{\text}
\def\f{\frac}

\def\em{\emptyset}

\def\bi{\binom}
\def\eq{\equiv}

\def\ls{\leqslant}
\def\gs{\geqslant}
\def\mo{\roman{mod}}

\def\Proof{\noindent{\it Proof}}

\def\Remark{\medskip\noindent{\it  Remark}}

\def\Ack{\medskip\noindent {\bf Acknowledgment}}
\hbox {Acta Arith. 156(2012), no.\,2, 123--141.}\bigskip \topmatter
\title On sums involving products of three binomial coefficients\endtitle
\author Zhi-Wei Sun\endauthor
\leftheadtext{Zhi-Wei Sun} \rightheadtext{Sums involving products of binomial coefficients}
\affil Department of Mathematics, Nanjing University
 \\Nanjing 210093, People's Republic of China
\\  zwsun\@nju.edu.cn\\ {\tt http://math.nju.edu.cn/$\sim$zwsun}
\endaffil

\abstract In this paper we mainly employ the Zeilberger algorithm to study congruences for sums of terms
involving products of three binomial coefficients.
Let $p>3$ be a prime. We prove that
$$\sum_{k=0}^{p-1}\f{\bi{2k}k^2\bi{2k}{k+d}}{64^k}\eq0\pmod{p^2}$$
for all $d\in\{0,\ldots,p-1\}$ with $d\eq(p+1)/2\pmod2$.
If $p\eq1\pmod4$ and $p=x^2+y^2$ with $x\eq1\pmod4$ and $y\eq0\pmod2$, then we show
$$\sum_{k=0}^{p-1}\f{\bi{2k}k^2\bi{2k}{k+1}}{(-8)^k}\eq2p-2x^2\ (\mo\ p^2)\ \t{and}\
\sum_{k=0}^{p-1}\f{\bi{2k}k\bi{2k}{k+1}^2}{(-8)^k}\eq-2p\ (\mo\ p^2)$$
by means of determining $x$ mod $p^2$ via
$$(-1)^{(p-1)/4}\,x\eq\sum_{k=0}^{(p-1)/2}\f{k+1}{8^k}\bi{2k}k^2
\eq\sum_{k=0}^{(p-1)/2}\f{2k+1}{(-16)^k}\bi{2k}k^2\pmod{p^2}.$$
We also solve the remaining open cases of Rodriguez-Villegas' conjectural congruences
on $$\sum_{k=0}^{p-1}\f{\bi{2k}k^2\bi{3k}k}{108^k},\ \ \sum_{k=0}^{p-1}\f{\bi{2k}k^2\bi{4k}{2k}}{256^k},
\ \ \sum_{k=0}^{p-1}\f{\bi{2k}{k}\bi{3k}k\bi{6k}{3k}}{12^{3k}}$$
modulo $p^2$.
\endabstract
\thanks 2010 {\it Mathematics Subject Classification}. Primary 11B65;
Secondary 05A10, 11A07, 11E25.
\newline\indent  Supported by the National Natural Science
Foundation (grant 11171140) of China and the PAPD of Jiangsu Higher
Education Institutions.
\endthanks
\endtopmatter
\document

\heading{1. Introduction}\endheading

Let $p$ be an odd prime.
It is known that (see, e.g., S. Ahlgren [A], L. van Hamme [vH] and T. Ishikawa [I])
$$\align&\sum_{k=0}^{(p-1)/2}(-1)^k\bi{-1/2}k^3
\\\eq&\cases 4x^2-2p\pmod{p^2}&\t{if}\ p=x^2+y^2\ (2\nmid x\ \&\ 2\mid y),
\\0\pmod{p^2}&\t{if}\ p\eq3\pmod{4}.\endcases
\endalign$$
Clearly,
$$\bi{-1/2}k=\f{\bi{2k}k}{(-4)^k}\quad\t{for all}\ k\in\N=\{0,1,2,\ldots\},$$
and
$$\bi{2k}k=\f{(2k)!}{(k!)^2}\eq0\pmod{p}\quad\t{for any}\ k=\f{p+1}2,\ldots,p-1.$$
After his determination of $\sum_{k=0}^{p-1}\bi{2k}k/m^k$ mod $p^2$
(where $m\in\Z$ and $m\not\eq0\pmod p$) in [Su1],
the author [Su2, Su3] posed some conjectures on $\sum_{k=0}^{p-1}\bi{2k}k^3/m^k$ mod $p^2$
with $m\in\{1,-8,16,-64,256,-512,4096\}$; for example, the author [Su2] conjectured that
$$\sum_{k=0}^{p-1}\bi{2k}k^3\eq\cases4x^2-2p\pmod{p^2}&\t{if}\ (\f p7)=1\ \&\ p=x^2+7y^2\ (x,y\in\Z),
\\0\pmod{p^2}&\t{if}\ (\f p7)=-1,\ \t{i.e.},\ p\eq 3,5,6\pmod7,\endcases\tag1.1$$
where $(-)$ denotes the Legendre symbol.
(It is known that if $(\f p7)=1$ then $p=x^2+7y^2$ for some $x,y\in\Z$; see, e.g., [C, p.\,31].)
Quite recently Z.-H. Sun [S2] made a certain
progress on those conjectures; in particular, he proved (1.1) in the case $(\f p7)=-1$
and confirmed the author's conjecture on $\sum_{k=0}^{p-1}\bi{2k}k^3/(-8)^k$ mod $p^2$.

 Let $p=2n+1$ be an odd prime. It is easy to see that for any $k=0,\ldots,n$ we have
 $$\bi{n+k}{2k}=\f{\prod_{j=1}^k(-(2j-1)^2)}{4^k(2k)!}\prod_{j=1}^k\l(1-\f{p^{2}}{(2j-1)^2}\r)
\eq\f{\bi{2k}k}{(-16)^k}\ (\mo\ p^2).
\tag1.2$$
Based on this observation Z.-H. Sun [S2] studied the polynomial
$$f_n(x)=\sum_{k=0}^n\bi{n+k}{2k}\bi{2k}k^2x^k$$
and found the key identity
$$f_n(x(x+1))=D_n(x)^2\tag1.3$$
in his approach to (1.1), where
$$D_n(x):=\sum_{k=0}^n\bi{n+k}{2k}\bi{2k}kx^k=\sum_{k=0}^n\bi nk\bi{n+k}kx^k.$$
Note that the numbers $D_n=D_n(1)\ (n\in\N)$ are the so-called central Delannoy numbers
and $P_n(x):=D_n((x-1)/2)$ is the Legendre polynomial of degree $n$.

Recall that the Catalan numbers are the integers defined by
$$C_n=\f1{n+1}\bi{2n}n=\bi{2n}n-\bi{2n}{n+1}\ \ \ (n\in\N)$$
while the Schr\"oder numbers are given by
$$S_n:=\sum_{k=0}^n\bi{n+k}{2k}C_k=\sum_{k=0}^n\bi nk\bi{n+k}k\f1{k+1}.$$
We define the Schr\"oder polynomial of degree $n$ by
$$S_n(x):=\sum_{k=0}^n\bi{n+k}{2k}C_kx^k.\tag1.4$$

For basic information about $D_n$ and $S_n$, the reader may consult
[CHV], [Sl],  [St, pp.\,178 and 185], and [Su4].

In combinatorics, Zeilberger's algorithm developed in [Z] (see also
Chapter 6 of [PWZ, pp.\,101-119]) is an algorithm which finds a
polynomial recurrence for a terminating hypergeometric sum. For
example, if we use {\tt Mathematica 7} and input {\tt
Zb[Binomial[n,k]${}^\wedge$3,\{k,0,n\},n,2]}, then we obtain the
following second-order recurrence for $S(n)=\sum_{k=0}^n\bi nk^3$:
$$-8(n+1)^2S(n)-(7n^2+21n+16)S(n+1)+(n+2)^2S(n+2)=0.$$

Via the Schr\"oder polynomials and the Zeilberger algorithm, we obtain
the following result.

\proclaim{Theorem 1.1} Let $p$ be an odd prime.

{\rm (i)} We have
$$\sum_{k=0}^{p-1}\f{\bi{2k}k^2\bi{2k}{k+d}}{64^k}\eq0\pmod{p^2}\tag1.5$$
for all $d\in\{0,1,\ldots,p-1\}$ with $d\eq(p+1)/2\pmod{2}$.

{\rm (ii)} If $p\eq3\pmod4$, then
$$\sum_{k=0}^{p-1}\f{\bi{2k}k^2\bi{2k}{k+1}}{64^k}
\eq(2p+2-2^{p-1})\bi{(p-1)/2}{(p+1)/4}^2\pmod{p^2}\tag1.6$$
\endproclaim

Now we state our second theorem the first part of which plays a key role in our proof of the second part.

\proclaim{Theorem 1.2} Let $p\eq1\pmod4$ be a prime and write $p=x^2+y^2$ with $x\eq1\pmod4$ and $y\eq0\pmod2$.

{\rm (i)} We can determine $x$ mod $p^2$ in the following way:
$$(-1)^{(p-1)/4}\,x\eq\sum_{k=0}^{(p-1)/2}\f{k+1}{8^k}\bi{2k}k^2
\eq\sum_{k=0}^{(p-1)/2}\f{2k+1}{(-16)^k}\bi{2k}k^2\pmod{p^2}.\tag1.7$$
Also,
$$\sum_{k=0}^{(p-1)/2}\f{\bi{2k}kC_k}{8^k}\eq-2\sum_{k=0}^{p-1}\f{k\bi{2k}k^2}{8^k}
\eq(-1)^{(p-1)/4}\l(2x-\f px\r)\pmod{p^2},\tag1.8$$
$$\aligned S_{(p-1)/2}\eq&\sum_{k=0}^{(p-1)/2}\f{\bi{2k}kC_k}{(-16)^k}\eq-8\sum_{k=0}^{(p-1)/2}\f{k\bi{2k}k^2}{(-16)^k}
\\\eq&(-1)^{(p-1)/4}\,2\l(2x-\f{p}x\r)\pmod{p^2},
\endaligned\tag1.9$$
$$\sum_{k=0}^{(p-1)/2}\f{k^2\bi{2k}k^2}{8^k}\eq(-1)^{(p-1)/4}\l(x-\f{3p}{4x}\r)\pmod{p^2},\tag1.10$$
and
$$\sum_{k=0}^{(p-1)/2}\f{k^2\bi{2k}k^2}{(-16)^{k}}\eq(-1)^{(p+3)/4}\f p{16x}\pmod{p^2}.\tag1.11$$

{\rm (ii)} We have
$$\sum_{k=0}^{p-1}\f{\bi{2k}k^2\bi{2k}{k+1}}{(-8)^k}\eq2p-2x^2\pmod{p^2}\tag1.12$$
and
$$\sum_{k=0}^{p-1}\f{\bi{2k}k\bi{2k}{k+1}^2}{(-8)^k}\eq-2p\pmod{p^2}.\tag1.13$$
\endproclaim
\Remark\ 1.1. Let $p$ be an odd prime. We conjecture that
$$\align&\sum_{k=0}^{p-1}\f{k+1}{8^k}\bi{2k}k^2+\sum_{k=0}^{(p-1)/2}\f{2k+1}{(-16)^k}\bi{2k}k^2
\\\eq&\cases 2(\f 2p)x\pmod{p^3}&\t{if}\ p=x^2+y^2\ (4\mid x-1\ \&\ 2\mid y),
\\0\pmod{p^2}&\t{if}\ p\eq3\pmod4.\endcases
\endalign$$
\medskip

Motivated by his study of Gaussian hypergeometric series and Calabi-Yau manifolds, in 2003 Rodriguez-Villegas [RV]
raised some conjectures on congruences. In particular, he conjectured that for any prime $p>3$ we have
$$\sum_{k=0}^{p-1}\f{\bi{2k}k^2\bi{3k}k}{108^k}\eq b(p)\pmod{p^2},
\ \ \ \ \sum_{k=0}^{p-1}\f{\bi{2k}k^2\bi{4k}{2k}}{256^k}\eq c(p)\pmod{p^2},\tag1.14$$
and
$$\sum_{k=0}^{p-1}\f{\bi{2k}k\bi{3k}k\bi{6k}{3k}}{12^{3k}}\eq\l(\f p3\r)a(p)\pmod{p^2},\tag1.15$$
where
$$\gather\sum_{n=1}^\infty a(n)q^n=q\prod_{n=1}^\infty(1-q^{4n})^6=\eta(4z)^6,
\\\sum_{n=1}^\infty b(n)q^n=q\prod_{n=1}^\infty(1-q^{6n})^3(1-q^{2n})^3=\eta^3(6z)\eta^3(2z),
\\\sum_{n=1}^\infty c(n)q^n=q\prod_{n=1}^\infty(1-q^n)^2(1-q^{2n})(1-q^{4n})(1-q^{8n})^2
=\eta^2(8z)\eta(4z)\eta(2z)\eta^2(z),
\endgather$$
and the Dedekind $\eta$-function is given by
$$\eta(z)=q^{1/24}\prod_{n=1}^\infty(1-q^n)\quad (\text{Im}(z)>0\ \t{and}\ q=e^{2\pi iz}).$$
In 1892 F. Klein and R. Fricke [KF] proved that (see also [SB])
$$a(p)=\cases 4x^2-2p&\t{if}\ p\eq1\pmod4\ \t{and}\ p=x^2+y^2\ (2\nmid x),
\\0&\t{if}\ p\eq3\pmod4.\endcases$$
By [SB] we also have
$$b(p)=\cases4x^2-2p&\t{if}\ p\eq1\pmod{3}\ \t{and}\ p=x^2+3y^2\ \t{with}\ x,y\in\Z,
\\0&\t{if}\ p\eq2\pmod3,\endcases$$
and
$$c(p)=\cases4x^2-2p&\t{if}\ (\f{-2}p)=1\ \t{and}\ p=x^2+2y^2\ \t{with}\ x,y\in\Z,
\\0&\t{if}\ (\f{-2}p)=-1,\ \t{i.e.,}\ p\eq5,7\pmod8.\endcases$$
Via an advanced approach involving the $p$-adic Gamma function and Gauss and Jacobi sums
(see K. Ono [O, Chapter 11] for an introduction to this method),
E. Mortenson [M] managed to provide a partial solution of (1.14) and (1.15), with the following congruences still open:
$$\gather\sum_{k=0}^{p-1}\f{\bi{2k}k^2\bi{3k}k}{108^k}\eq b(p)=0\pmod{p^2}\quad\t{if}\ p\eq5\pmod6,\tag1.16
\\\sum_{k=0}^{p-1}\f{\bi{2k}k^2\bi{4k}{2k}}{256^k}\eq c(p)\pmod{p^2}\quad\t{if}\ p\eq3\pmod4,\tag1.17
\\\sum_{k=0}^{p-1}\f{\bi{2k}k\bi{3k}{k}\bi{6k}{3k}}{12^{3k}}\eq -a(p)\pmod{p^2}\quad\t{if}\ p\eq5\pmod6.\tag1.18
\endgather$$
Concerning (1.16)-(1.18), Mortenson's approach [M] only allowed him
to show that for each of them the squares of both sides of the congruence
 are congruent modulo $p^2$.

 Our following theorem confirms (1.16)-(1.18) and hence completes the proof of (1.14) and (1.15).
So far, all conjectures of Rodriguez-Villegas [RV] involving at most three products of binomial coefficients
have been proved!

 \proclaim{Theorem 1.3} Let $p>3$ be a prime.

 {\rm (i)} Given $d\in\{0,\ldots,p-1\}$, we have
$$\gather\sum_{k=0}^{p-1}\f{\bi{2k}{k+d}\bi{2k}k\bi{3k}{k}}{108^k}\eq0\pmod{p^2}\quad\t{if}\ d\eq \f{1+(\f p3)}2\pmod{2},\tag1.19
\\\sum_{k=0}^{p-1}\f{\bi{2k}{k+d}\bi{2k}k\bi{4k}{2k}}{256^k}\eq0\pmod{p^2}\quad\t{if}\ d\eq\f{1+(\f{-2}p)}2\pmod{2},\tag1.20
\\\sum_{k=0}^{p-1}\f{\bi{2k}{k+d}\bi{3k}k\bi{6k}{3k}}{12^{3k}}\eq0\pmod{p^2}\quad\t{if}\ d\eq\f{1+(\f{-1}p)}2\pmod{2}.\tag1.21
\endgather$$

{\rm (ii)} If $p\eq3\pmod8$ and $p=x^2+2y^2$ with $x,y\in\Z$, then
$$\sum_{k=0}^{p-1}\f{\bi{2k}k^2\bi{4k}{2k}}{256^k}\eq4x^2-2p\pmod{p^2}.\tag1.22$$

{\rm (iii)} If $p\eq5\pmod{12}$ and $p=x^2+y^2$ with $2\nmid x$ and $2\mid y$, then
$$\sum_{k=0}^{p-1}\f{\bi{2k}k\bi{3k}k\bi{6k}{3k}}{12^{3k}}\eq2p-4x^2\pmod{p^2}.\tag1.23$$
\endproclaim

In the case $d=1$, Theorem 1.3(i) yields the following new result. (Note that $\bi{2k}k\bi{3k}{k+1}=2\bi{2k}{k+1}\bi{3k}k$
for any $k\in\N$.)

\proclaim{Corollary 1.1} Let $p>3$ be a prime. Then
$$\gather\sum_{k=0}^{p-1}\f{\bi{2k}k^2\bi{3k}{k+1}}{108^k}\eq0\pmod{p^2}\quad\t{if}\ p\eq1\pmod{3},\tag1.24
\\\sum_{k=0}^{p-1}\f{\bi{4k}{2k}\bi{2k}k\bi{2k}{k+1}}{256^k}\eq0\pmod{p^2}\quad\t{if}\ p\eq1,3\pmod{8},\tag1.25
\\\sum_{k=0}^{p-1}\f{\bi{6k}{3k}\bi{3k}k\bi{2k}{k+1}}{12^{3k}}\eq0\pmod{p^2}\quad\t{if}\ p\eq1\pmod{4}.\tag1.26
\endgather$$
\endproclaim

We will prove Theorems 1.1-1.3 in Sections 2-4 respectively.

\heading{2. Proof of Theorem 1.1}\endheading

\proclaim{Lemma 2.1} For any positive integer $n$ we have
$$\sum_{k=1}^n\bi{n+k}{2k}\bi{2k}k\bi{2k}{k+1}x^{k-1}(x+1)^{k+1}=n(n+1)S_n(x)^2.\tag2.1$$
\endproclaim
\Proof. Observe that
$$S_n(x)^2=\sum_{k=0}^n\bi{n+k}{2k}C_kx^k\sum_{l=0}^n\bi{n+l}{2l}C_lx^l=\sum_{m=0}^{2n}a_m(n)x^m,$$
where
$$a_m(n):=\sum_{k=0}^m\bi{n+k}{2k}C_k\bi{n+m-k}{2m-2k}C_{m-k}.$$
Also, the coefficient of $x^m$ on the left-hand side of (2.1) coincides with
$$\align b_m(n):=&\sum_{k=1}^{m+1}\bi{n+k}{2k}\bi{2k}k\bi{2k}{k+1}\bi{k+1}{m+1-k}
\\=&\sum_{k=0}^m\bi{n+k+1}{2k+2}\bi{2k+2}{k+1}\bi{2k+2}k\bi{k+2}{m-k}.
\endalign$$
Thus, for the validity of (2.1) it suffices to show that $b_m(n)=n(n+1)a_m(n)$ for all $m=0,1,\ldots$.
Obviously, $a_0(n)=1$ and $b_0(n)=n(n+1)$. Also, $a_1(n)=n(n+1)$ and $b_1(n)=n^2(n+1)^2$.
By the Zeilberger algorithm via {\tt Mathematica 7} we find that both $u_m=a_m(n)$ and $u_m=b_m(n)$
satisfy the following recursion:
$$\align&(m+2)(m+3)(m+4)u_{m+2}
\\=&2(2mn^2+5n^2+2mn+5n-m^3-6m^2-11m-6)u_{m+1}
\\&-(m+1)(m-2n)(m+2n+2)u_m.
\endalign$$
So $b_m(n)=n(n+1)a_m(n)$ for all $m\in\N$. This proves (2.1). \qed

\medskip
\noindent{\it Proof of Theorem} 1.1.
 We first determine $\sum_{k=0}^{p-1}\bi{2k}k^2\bi{2k}{k+1}/64^k$ mod $p^2$ via Lemma 2.1,
which actually led the author to the study of (1.5).

Recall the following combinatorial identity (cf. [Su2, (4.3)]):
$$\sum_{k=0}^n\bi{n+k}{2k}\f{C_k}{(-2)^k}=\cases(-1)^{(n-1)/2}C_{(n-1)/2}/2^n&\t{if}\ 2\nmid n,
\\0&\t{if}\ 2\mid n.\endcases\tag2.2$$

Set $n=(p-1)/2$. Applying (2.1) with $x=-1/2$ we get
$$\sum_{k=1}^n\bi{n+k}{2k}\bi{2k}k\bi{2k}{k+1}\f1{(-2)^{k-1}2^{k+1}}=n(n+1)S_n\l(-\f12\r)^2.$$
Thus, with the helps from (1.2) and (2.2), we have
$$\align\sum_{k=0}^{p-1}\f{\bi{2k}k^2\bi{2k}{k+1}}{64^k}
\eq&\sum_{k=1}^n\bi{n+k}{2k}\bi{2k}k\bi{2k}{k+1}\f1{(-4)^k}
\\=&-n(n+1)S_n\l(-\f12\r)^2
\\\eq&\cases0\pmod{p^2}&\t{if}\ p\eq1\pmod{4}
\\C_{(n-1)/2}^2/2^{2n+2}\pmod{p^2}&\t{if}\ p\eq3\pmod{4}.\endcases
\endalign$$
Therefore (1.5) with $d=1$ holds if $p\eq1\pmod4$.
In the case $p\eq3\pmod4$, clearly
$$\align\f{C_{(n-1)/2}^2}{2^{2n+2}}=&\f{\l(\bi{(p-1)/2}{(p+1)/4}\f2{p-1}\r)^2}{4\times 2^{p-1}}
\\\eq&\f1{(1-2p)(1+p\,q_p(2))}\bi{(p-1)/2}{(p+1)/4}^2
\\\eq&(1+2p-p\,q_p(2))\bi{(p-1)/2}{(p+1)/4}^2\pmod{p^2}
\endalign$$
where $q_p(2)=(2^{p-1}-1)/p$, and hence (1.6) holds.

For $d=0,1,2,\ldots$ set
$$u_d=\sum_{k=0}^{p-1}\f{\bi{2k}k^2\bi{2k}{k+d}}{64^k}=\sum_{d\ls k<p}\f{\bi{2k}k^2\bi{2k}{k+d}}{64^k}.$$
By the Zeilberger algorithm we find the recursion
$$(2d+1)^2u_d-(2d+3)^2u_{d+2}=\f{(2p-1)^2(d+1)}{64^{p-1}p}\bi{2p}{p+d+1}\bi{2p-2}{p-1}^2.$$
Note that
$$\bi{2p-2}{p-1}=pC_{p-1}\eq0\pmod{p}.$$
If $0\ls d<p-2$, then
$$\bi{2p}{p+d+1}=\f{2p}{p+d+1}\bi{2p-1}{p+d}\eq0\pmod{p}$$
and hence
$$(2d+1)^2u_d\eq (2d+3)^2u_{d+2}\pmod{p^2}.$$
For $d\in\{0,\ldots,p-3\}$ with $d\eq(p+1)/2\pmod2$, clearly $p\not=2d+1<2p$ and hence
$$u_{d+2}\eq0\pmod{p^2}\ \Longrightarrow\ u_d\eq0\pmod{p^2}.$$
If $d\in\{p-1,p-2\}$ and $d\eq(p+1)/2\pmod2$, then $d\gs(p+1)/2$ and hence $u_d\eq0\pmod{p^2}$.
So (1.5) holds for all $d\in\{0,\ldots,p-1\}$ with $d\eq(p+1)/2\pmod2$.

So far we have completed the proof of Theorem 1.1.
\qed

\heading{3. Proof of Theorem 1.2}\endheading

\proclaim{Lemma 3.1} For any $n\in\N$ we have
$$\sum_{k=0}^n\bi{2k}k^3\bi k{n-k}(-16)^{n-k}=\sum_{k=0}^n\bi{2k}k^2\bi{2(n-k)}{n-k}^2.\tag3.1$$
\endproclaim
\Proof. For $n=0,1$, both sides of (3.1) take the values $1$ and $8$ respectively.
Let $u_n$ denote the left-hand side of (3.1) or the right-hand side of (3.1).
Applying the Zeilberger algorithm via {\tt Mathematica 7}, we obtain the recursion
$$(n+2)^3u_{n+2}=8(2n+3)(2n^2+6n+5)u_{n+1}-256(n+1)^3u_n\ (n\in\N).$$
So, by induction (3.1) holds for all $n=0,1,2,\ldots$. \qed

\proclaim{Lemma 3.2} Let $p$ be an odd prime. Then
$$\align&\sum_{n=0}^{p-1}\f{n+1}{8^n}\sum_{k=0}^n\bi{2k}k^2\bi{2(n-k)}{n-k}^2
\\\eq&\sum_{n=0}^{p-1}\f{2n+1}{(-16)^n}\sum_{k=0}^n\bi{2k}k^2\bi{2(n-k)}{n-k}^2
\\\eq&p\l(\f{-1}p\r)\pmod{p^3}.
\endalign$$
\endproclaim
\Proof. In view of Lemma 3.1, we have
$$\align&\sum_{n=0}^{p-1}\f{n+1}{8^n}\sum_{k=0}^n\bi{2k}k^2\bi{2(n-k)}{n-k}^2
\\=&\sum_{n=0}^{p-1}\f{n+1}{8^n}\sum_{k=0}^n\bi{2k}k^3\bi k{n-k}(-16)^{n-k}
\\=&\sum_{k=0}^{p-1}\f{\bi{2k}k^3}{8^k}\sum_{j=0}^{p-1-k}(k+j+1)\bi kj\f{(-16)^j}{8^j}
\\\eq&\sum_{k=0}^{(p-1)/2}\f{\bi{2k}k^3}{8^k}\((k+1)\sum_{j=0}^k\bi kj(-2)^j-2k\sum_{j=1}^{k}\bi {k-1}{j-1}(-2)^{j-1}\)
\\=&\sum_{k=0}^{(p-1)/2}\f{\bi{2k}k^3}{8^k}\l((k+1)(-1)^k-2k(-1)^{k-1}\r)
\\\eq&\sum_{k=0}^{p-1}\f{3k+1}{(-8)^k}\bi{2k}k^3\pmod{p^3}.
\endalign$$
In [Su3] the author conjectured that
$$\sum_{k=0}^{p-1}\f{3k+1}{(-8)^k}\bi{2k}k^3\eq p\l(\f{-1}p\r)+p^3E_{p-3}\pmod{p^4}$$
provided $p>3$, where $E_0,E_1,E_2,\ldots$ are Euler numbers given by
$$E_0=1\ \ \t{and}\ \ \sum^n\Sb k=0\\2\mid k\endSb \bi nk E_{n-k}=0\ \ (n=1,2,3,\ldots).$$
The last congruence is still open but [GZ] confirmed that
$$\sum_{k=0}^{p-1}\f{3k+1}{(-8)^k}\bi{2k}k^3\eq p\l(\f{-1}p\r)\pmod{p^3}.$$
So we have
$$\sum_{n=0}^{p-1}\f{n+1}{8^n}\sum_{k=0}^n\bi{2k}k^2\bi{2(n-k)}{n-k}^2\eq p\l(\f{-1}p\r)\pmod{p^3}.$$
Similarly,
$$\align&\sum_{n=0}^{p-1}\f{2n+1}{(-16)^n}\sum_{k=0}^n\bi{2k}k^2\bi{2(n-k)}{n-k}^2
\\=&\sum_{n=0}^{p-1}\f{2n+1}{(-16)^n}\sum_{k=0}^n\bi{2k}k^3\bi k{n-k}(-16)^{n-k}
\\\eq&\sum_{k=0}^{(p-1)/2}\f{\bi{2k}k^3}{(-16)^k}\((2k+1)\sum_{j=0}^k\bi kj+2k\sum_{j=1}^k\bi{k-1}{j-1}\)
\\\eq&\sum_{k=0}^{p-1}\f{3k+1}{(-8)^k}\bi{2k}k^3\eq p\l(\f{-1}p\r)\pmod{p^3}.
\endalign$$
This concludes the proof. \qed

\proclaim{Lemma 3.3} Let $p$ be an odd prime. Then
$$\align2\sum_{k=0}^{(p-1)/2}\f{k\bi{2k}k^2}{8^k}+\sum_{k=0}^{(p-1)/2}\f{\bi{2k}kC_k}{8^k}
\eq& 2p^2\l(\f 2p\r)\pmod{p^3},
\\8\sum_{k=0}^{(p-1)/2}\f{k\bi{2k}k^2}{(-16)^k}+\sum_{k=0}^{(p-1)/2}\f{\bi{2k}kC_k}{(-16)^k}
\eq&2p^2\l(\f {-1}p\r)\pmod{p^3},
\\\sum_{k=0}^{(p-1)/2}(2k^2+4k+1)\f{\bi{2k}k^2}{8^k}\eq&p^2\l(\f 2p\r)\pmod{p^3},
\\\sum_{k=0}^{(p-1)/2}(8k^2+4k+1)\f{\bi{2k}k^2}{(-16)^k}\eq&p^2\l(\f {-1}p\r)\pmod{p^3}.
\endalign$$
\endproclaim
\Proof. By induction, for every $n=0,1,2,\ldots$ we have
$$\align\sum_{k=0}^{n}\l(2k+\f1{k+1}\r)\f{\bi{2k}k^2}{8^k}=&\f{(2n+1)^2}{(n+1)8^{n}}\bi{2n}{n}^2,
\\\sum_{k=0}^{n}\l(8k+\f1{k+1}\r)\f{\bi{2k}k^2}{(-16)^k}=&\f{(2n+1)^2}{(n+1)(-16)^{n}}\bi{2n}{n}^2,
\\\sum_{k=0}^n(2k^2+4k+1)\f{\bi{2k}k^2}{8^k}=&\f{(2n+1)^2}{8^n}\bi{2n}n^2,
\\\sum_{k=0}^n(8k^2+4k+1)\f{\bi{2k}k^2}{(-16)^k}=&\f{(2n+1)^2}{(-16)^n}\bi{2n}n^2.
\endalign$$
Applying these identities with $n=(p-1)/2$ we immediately get the desired congruences. \qed
\medskip

Let $p\eq1\pmod 4$ be a prime and write $p=x^2+y^2$ with $x\eq1\pmod4$ and $y\eq0\pmod2$.
In 1828 Gauss showed the congruence $\bi{(p-1)/2}{(p-1)/4}\eq2x\pmod p$.
In 1986, S. Chowla, B. Dwork and R. J. Evans [CDE] used Gauss and Jacobi sums to prove that
$$\bi{(p-1)/2}{(p-1)/4}\eq\f{2^{p-1}+1}2\l(2x-\f p{2x}\r)\pmod {p^2},\tag3.2$$
which was first conjectured by F. Beukers. (See also [BEW, Chapter 9] and [HW] for further related results.)
In 2009, the author (see [Su2]) conjectured that
$$\sum_{k=0}^{(p-1)/2}\f{\bi{2k}k^2}{8^k}\eq\sum_{k=0}^{(p-1)/2}\f{\bi{2k}k^2}{(-16)^k}
\eq(-1)^{(p-1)/4}\l(2x-\f p{2x}\r)\pmod{p^2},\tag3.3$$
and this was confirmed by Z.-H. Sun [S1] via (3.2) and the Legendre polynomials.

\medskip
\noindent{\it Proof of Theorem} 1.2(i). By (1.2),
$$S_{(p-1)/2}\eq\sum_{k=0}^{(p-1)/2}\f{\bi{2k}kC_k}{(-16)^k}\pmod{p^2}.$$
In view of this and Lemma 3.3 and (3.3), it suffices to show (1.7).

As $p\mid\bi{2k}k$ for all $k=(p+1)/2,\ldots,p-1$, we have
$$\align&\sum_{n=0}^{p-1}\f{n+1}{8^n}\sum_{k=0}^n\bi{2k}k^2\bi{2(n-k)}{n-k}^2
\\=&\sum_{k=0}^{p-1}\f{\bi{2k}k^2}{8^k}\sum_{n=k}^{p-1}\f{n+1}{8^{n-k}}\bi{2(n-k)}{n-k}^2
\\=&\sum_{k=0}^{p-1}\f{\bi{2k}k^2}{8^k}\sum_{j=0}^{p-1-k}\f{k+j+1}{8^{j}}\bi{2j}{j}^2
\\\eq&\sum_{k=0}^{(p-1)/2}\f{\bi{2k}k^2}{8^k}\sum_{j=0}^{(p-1)/2}\f{(k+1)+(j+1)-1}{8^j}\bi{2j}j^2
\\=&2\sum_{k=0}^{(p-1)/2}\f{\bi{2k}k^2}{8^k}\sum_{j=0}^{(p-1)/2}\f{(j+1)\bi{2j}j^2}{8^j}
-\(\sum_{k=0}^{p-1}\f{\bi{2k}k^2}{8^k}\)^2\pmod{p^2}.
\endalign$$
Similarly,
$$\align&\sum_{n=0}^{p-1}\f{2n+1}{(-16)^n}\sum_{k=0}^n\bi{2k}k^2\bi{2(n-k)}{n-k}^2
\\\eq&2\sum_{k=0}^{(p-1)/2}\f{\bi{2k}k^2}{(-16)^k}\sum_{j=0}^{(p-1)/2}\f{(2j+1)\bi{2j}j^2}{(-16)^j}
-\(\sum_{k=0}^{p-1}\f{\bi{2k}k^2}{(-16)^k}\)^2\pmod{p^2}.
\endalign$$
Combining these with Lemma 3.2 and (3.3), we immediately obtain (1.7). \qed

\proclaim{Lemma 3.4}  Let $p\eq1\pmod4$ be a prime. Write $p=x^2+y^2$ with $x\eq1\pmod4$ and $y\eq0\pmod2$. Then
$$D_{(p-1)/2}\eq(-1)^{(p-1)/4}\l(2x-\f p{2x}\r)\pmod{p^2}.\tag3.4$$
\endproclaim
\Proof. By (1.2),
$$D_{(p-1)/2}\eq\sum_{k=0}^{(p-1)/2}\f{\bi{2k}k^2}{(-16)^k}\pmod{p^2}.$$
So (3.4) follows from (3.3). \qed
\Remark\ 3.1. If $p$ is a prime with $p\eq3\pmod4$, then $n=(p-1)/2$ is odd and hence
$$\align D_n\eq&\sum_{k=0}^n(-1)^k\f{\bi{2k}k^2}{16^k}
=\sum_{k=0}^n(-1)^k\bi{-1/2}k^2
\\\eq&\sum_{k=0}^n(-1)^k\bi{n}k^2=\sum_{k=0}^n(-1)^{n-k}\bi{n}{k}^2=0\pmod p.
\endalign$$

The following result was conjectured by the author [Su2] and confirmed by Z.-H. Sun [S2].
\proclaim{Lemma 3.5}  Let $p$ be an odd prime. Then
$$\sum_{k=0}^{p-1}\f{\bi{2k}k^3}{(-8)^k}\eq\cases 4x^2-2p\ (\mo\ p^2)&\t{if}\ 4\mid p-1\ \&\ p=x^2+y^2\ (2\nmid x),
\\0\ (\mo\ p^2)&\t{if}\ p\eq3\ (\mo\ 4).\endcases\tag3.5$$
\endproclaim
\Remark\ 3.2. Fix an odd prime $p=2n+1$. By (1.2) and (1.3) we have
$$\sum_{k=0}^{p-1}\f{\bi{2k}k^3}{(-8)^k}\eq\sum_{k=0}^n\bi{n+k}{2k}\bi{2k}k^22^k=D_n^2\pmod{p^2}.$$
 Hence (3.5) follows from Lemma 3.4 and Remark 3.1.

\proclaim{Lemma 3.6} For any positive integer $n$ we have
$$\sum_{k=0}^n\bi{n+k}{2k}\bi{2k}k^2\f{2k+1}{(k+1)^2}x^k(x+1)^{k+1}=\f{S_n(x)}2(D_{n-1}(x)+D_{n+1}(x)).\tag3.6$$
\endproclaim
\Proof. Note that
$$S_n(x)(D_{n-1}(x)+D_{n+1}(x))=\sum_{m=0}^{2n+1}c_m(n)x^m$$
where
$$\align &c_m(n)=\sum_{k=0}^m\bi{n+k}{2k}C_k\bi{2m-2k}{m-k}\(\bi{n-1+m-k}{2m-2k}+\bi{n+1+m-k}{2m-2k}\)
\\=&2\sum_{k=0}^m\bi{n+k}{2k}C_k\bi{n+m-k}{2m-2k}\bi{2m-2k}{m-k}\f{(m+n-k)^2-n(2m-2k-1)}{(m+n-k)(n-m+k+1)}.
\endalign$$
By the Zeilberger algorithm we find that $u_m=c_m(n)/2$ satisfies the recursion
$$\aligned &(m+2)(m+3)^2(m^2+5m+6+4n(n+1))u_{m+2}+2P(m,n)u_{m+1}
\\=&(m+2)((2n+1)^2-m^2)(m^2+7m+12+4n(n+1))u_m
\endaligned\tag3.7$$
where $P(m,n)$ denotes the polynomial
$$\align &m^5+11m^4+45m^3+83m^2+64m+12+20n^4-40n^3-58n^2-38n
\\&-25mn+m^2n+2m^3n-33mn^2+m^2n^2+2m^3n^2-16mn^3-8mn^4.
\endalign$$
Clearly the coefficient of $x^m$ on the left-hand side of (3.6) coincides with
$$d_m(n)=\sum_{k=0}^m\bi{n+k}{2k}\bi{2k}k^2\bi{k+1}{m-k}\f{2k+1}{(k+1)^2}.$$
By the Zeilberger algorithm $u_m=d_m(n)$ also satisfies the recursion (3.7).
Thus we have $d_m(n)=c_m(n)$ by induction on $m$. So (3.6) holds. \qed

\medskip
\noindent{\it Proof of Theorem} 1.2(ii). Write $p=2n+1$. By (2.1),
$$\sum_{k=0}^n\bi{n+k}{2k}\bi{2k}k\bi{2k}{k+1}2^k=\f{n(n+1)}2S_n^2.$$
Thus, by (1.2) and (1.9) we have
$$\align\sum_{k=0}^{p-1}\f{\bi{2k}k^2\bi{2k}{k+1}}{(-8)^k}
\eq&\sum_{k=0}^n\bi{n+k}{2k}\bi{2k}k\bi{2k}{k+1}2^k
\\\eq&\f{p^2-1}{8}4(4x^2-4p)\pmod{p^2}
\endalign$$
and hence (1.12) holds.

Now we consider (1.13). Observe that
$$\bi{2k}{k+1}^2=\l(1-\f{2k+1}{(k+1)^2}\r)\bi{2k}k^2\quad\t{for}\ k=0,1,2,\ldots,$$
and
$$\bi{2(p-1)}{p-1}\bi{2(p-1)}{(p-1)+1}^2=\f p{2p-1}\bi{2p-1}{p-1}\bi{2p-2}{p-2}^2\eq-p\pmod{p^2}.$$
Thus we have
$$\sum_{k=0}^{p-1}\f{\bi{2k}k\bi{2k}{k+1}^2}{(-8)^k}\eq-p+\sum_{k=0}^n\f{\bi{2k}k^3}{(-8)^k}
-\sum_{k=0}^n\f{(2k+1)\bi{2k}k^3}{(k+1)^2(-8)^k}\pmod{p^2}.\tag3.8$$

By (1.2) and (3.6) with $x=1$,
$$\align \sum_{k=0}^n\f{(2k+1)\bi{2k}k^3}{(k+1)^2(-8)^k}\eq&\sum_{k=0}^n\bi{n+k}{2k}\bi{2k}k^2\f{(2k+1)2^k}{(k+1)^2}
\\=&\f{S_n}4(D_{n-1}+D_{n+1})\pmod{p^2}.
\endalign$$
It is known (cf. [Sl] and [St, p.\,191]) that
$$(n+1)D_{n+1}=3(2n+1)D_n-nD_{n-1}\quad\t{and}\quad D_{n+1}-3D_n=2nS_n.$$
Thus
$$\align n(D_{n-1}+D_{n+1})=&3(2n+1)D_n-D_{n+1}
\\=&3(2n+1)D_n-(3D_n+2nS_n)=2n(3D_n-S_n)
\endalign$$
and hence
$$\sum_{k=0}^n\f{(2k+1)\bi{2k}k^3}{(k+1)^2(-8)^k}\eq\f{S_n}2(3D_n-S_n)\pmod{p^2}.$$

With the help of (1.9) and (3.4), we have
$$\f{S_n}2(3D_n-S_n)\eq\l(2x-\f px\r)\l(3\l(2x-\f p{2x}\r)-\l(4x-\f{2p}x\r)\r)\pmod{p^2}$$
and hence
$$\sum_{k=0}^n\f{(2k+1)\bi{2k}k^3}{(k+1)^2(-8)^k}\eq4x^2-p\pmod{p^2}.$$
Combining this with (3.5) and (3.8), we immediately obtain (1.13). \qed

\heading{4. Proof of Theorem 1.3}\endheading

\proclaim{Lemma 4.1} Let $p$ be an odd prime. Then, for any $p$-adic integer $x\not\eq0,-1\ (\mo\ p)$ we have
$$\sum_{k=0}^{p-1}\bi{2k}k^3\l(\f {-x}{64}\r)^k
\eq\l(\f{x+1}p\r)\sum_{k=0}^{p-1}\bi{2k}k^2\bi{4k}{2k}\l(\f x{64(x+1)^2}\r)^k
\ (\mo\ p).\tag4.1$$
\endproclaim
\Proof. Taking $n=(p-1)/2$ in the following identity of
MacMahon (see, e.g., [G, (6.7)])
$$\sum_{k=0}^n\bi nk^3x^k=\sum_{k=0}^n\bi {n+k}{2k}\bi{2k}k\bi{n-k}kx^k(1+x)^{n-2k}$$
and noting (1.2) and the basic facts
$$\bi nk\eq\bi{-1/2}k=\f{\bi{2k}k}{(-4)^k}\pmod p$$
and
$$\bi {n-k}k\eq\bi{-1/2-k}k=\f{\bi{4k}{2k}}{(-4)^k}\pmod p,$$
we immediately get (4.1). \qed

\medskip
\noindent{\it Proof of Theorem 1.3}. (i) For $d=0,1,2,\ldots$, we
define
$$f(d)=\sum_{k=0}^{p-1}\f{\bi{2k}{k+d}\bi{2k}k\bi{3k}k}{108^k},\ \ \
g(d)=\sum_{k=0}^{p-1}\f{\bi{2k}{k+d}\bi{2k}k\bi{4k}{2k}}{256^k},$$
and
$$h(d)=\sum_{k=0}^{p-1}\f{\bi{2k}{k+d}\bi{3k}k\bi{6k}{3k}}{12^{3k}}.$$
By the Zeilberger algorithm, we find the recursive relations:
$$\aligned &(3d+1)(3d+2)f(d)-(3d+4)(3d+5)f(d+2)
\\=&\f{(3p-1)(3p-2)(d+1)}{108^{p-1}p}\bi{2p}{p+d+1}\bi{2p-2}{p-1}\bi{3p-3}{p-1},
\endaligned\tag4.2$$
$$\aligned &(4d+1)(4d+3)g(d)-(4d+5)(4d+7)g(d+2)
\\=&\f{(4p-1)(4p-3)(d+1)}{256^{p-1}p}\bi{2p}{p+d+1}\bi{2p-2}{p-1}\bi{4p-4}{2p-2},
\endaligned\tag4.3$$
and
$$\aligned &(6d+1)(6d+5)h(d)-(6d+7)(6d+11)h(d+2)
\\=&\f{(6p-1)(6p-5)(d+1)}{1728^{p-1}p}\bi{2p}{p+d+1}\bi{3p-3}{p-1}\bi{6p-6}{3p-3}.
\endaligned\tag4.4$$
Recall that $\bi{2p-2}{p-1}=pC_{p-1}\eq0\pmod{p}$. Also,
$$\align (3p-2)\bi{3p-3}{p-1}=&p\bi{3p-2}p\eq0\pmod p,
\\(4p-3)\bi{4p-4}{2p-2}=&p\bi{4p-2}{2p}\eq0\pmod{p},
\\(6p-5)\bi{6p-6}{3p-3}=&\f{3p(3p-1)(3p-2)}{(6p-3)(6p-4)}\bi{6p-3}{3p}\eq0\pmod p.
\endalign$$
If $0\ls d<p-1$, then
$$\bi{2p}{p+d+1}=\bi{2p}{p-1-d}\eq0\pmod p.$$
So, by (4.2)-(4.4), for any $d\in\{0,\ldots,p-1\}$ we have
$$\align(3d+1)(3d+2)f(d)\eq&(3d+4)(3d+5)f(d+2)\pmod{p^2},\tag4.5
\\(4d+1)(4d+3)g(d)\eq&(4d+5)(4d+7)g(d+2)\pmod{p^2},\tag4.6
\\(6d+1)(6d+5)h(d)\eq&(6d+7)(6d+11)h(d+2)\pmod{p^2}.\tag4.7
\endalign$$

Fix $0\ls d\ls p-1$. If $d\eq(1+(\f p3))/2\pmod2$, then it is easy
to verify that $\{3d+1,3d+2\}\cap\{p,2p\}=\em$, hence
$(3d+1)(3d+2)\not\eq0\pmod{p}$ and thus by (4.5) we have
$$f(d+2)\eq0\ (\mo\ p^2)\ \Longrightarrow\ f(d)\eq0\ (\mo\ p^2).$$
If $d\eq(1+(\f {-2}p))/2\pmod2$, then
$\{4d+1,4d+3\}\cap\{p,3p\}=\em$, hence
$(4d+1)(4d+3)\not\eq0\pmod{p}$ and thus by (4.6) we have
$$g(d+2)\eq0\ (\mo\ p^2)\ \Longrightarrow\ g(d)\eq0\ (\mo\ p^2).$$
If $d\eq(1+(\f {-1}p))/2\pmod2$, then
$\{6d+1,6d+3\}\cap\{p,3p,5p\}=\em$, hence
$(6d+1)(6d+3)\not\eq0\pmod{p}$ and thus (4.7) yields
$$h(d+2)\eq0\ (\mo\ p^2)\ \Longrightarrow\ h(d)\eq0\ (\mo\ p^2).$$

Since
$$f(p)=f(p+1)=g(p)=g(p+1)=h(p)=h(p+1)=0,$$
by the last paragraph, for every $d=p+1,p,\ldots,0$ we have the desired
(1.19)-(1.21).

(ii) Assume that $p\eq3\pmod8$ and $p=x^2+2y^2$ with $x,y\in\Z$. Since $4x^2\not\eq0\pmod p$ and
Mortenson [M] already proved that the squares of both sides of (1.22) are congruent modulo $p^2$,
(1.22) is reduced to its mod $p$ form. Applying (4.1) with $x=1$ we get
$$\sum_{k=0}^{p-1}\f{\bi{2k}k^3}{(-64)^k}\eq
\l(\f2p\r)\sum_{k=0}^{p-1}\f{\bi{2k}k^2\bi{4k}{2k}}{256^k}\pmod p.$$
By [A, Theorem 5(3)], we have
$$\l(\f{-1}p\r)\sum_{k=0}^n\bi nk^2\bi{n+k}k(-1)^k\eq 4x^2-2p\pmod{p},$$
where $n=(p-1)/2$.  For $k=0,\ldots,n$ clearly
$$\align\bi nk^2\bi{n+k}k(-1)^k=&\bi{(p-1)/2}k^2\bi{-(p+1)/2}k
\\\eq&\bi{-1/2}k^3=\f{\bi{2k}k^3}{(-64)^k}\pmod{p},
\endalign$$
therefore $$\sum_{k=0}^{p-1}\f{\bi{2k}k^3}{(-64)^k}\eq \l(\f{-1}p\r)(4x^2-2p)\pmod {p}$$
and hence (1.22) follows.

(iii) Finally we suppose $p\eq5\pmod{12}$ and write $p=x^2+y^2$ with $x$ odd and $y$ even. Once again
it suffices to show the mod $p$ form of (1.23) in view of Mortenson's work [M].
As Z.-H. Sun observed,
$$\bi{(p-5)/6+k}{2k}\bi{2k}k\eq\bi{k-5/6}{2k}\bi{2k}k=\f{\bi{3k}k\bi{6k}{3k}}{(-432)^k}\pmod p$$
for all $k=0,1,2,\ldots$. If $p/6<k<p/3$ then $p\mid\bi{6k}{3k}$; if $p/3<k<p/2$ then $p\mid\bi{3k}k$;
if $p/2<k<p$ then $p\mid\bi{2k}k$.
Thus
$$\align\sum_{k=0}^{p-1}\f{\bi{2k}k\bi{3k}k\bi{6k}{3k}}{12^{3k}}
\eq&\sum_{k=0}^{(p-5)/6}\bi{(p-5)/6+k}{2k}\bi{2k}k^2\l(-\f14\r)^k
\\=&D_{2n}\l(-\f12\r)^2\pmod{p}\quad(\t{by (1.3)}),
\endalign$$
where $n=(p-5)/12$. Note that
$$D_{2n}\l(-\f12\r)=\f1{(-4)^n}\bi{2n}n$$
by [G, (3.133) and (3.135)], and
$$\bi{(p-1)/2}{(p-1)/4}\eq12(-432)^n\bi{2n}n\pmod{p}$$
by P. Morton [Mo].
Therefore
$$D_{2n}\l(-\f12\r)^2=\f1{16^n}\bi{2n}n^2\eq\f{\bi{(p-1)/2}{(p-1)/4}^2}{12^{6n+2}}\eq\l(\f{12}p\r)\bi{(p-1)/2}{(p-1)/4}^2
\pmod{p}.$$
Thus, by applying Gauss' congruence $\bi{(p-1)/2}{(p-1)/4}\eq2x\pmod{p}$ (cf. [BEW, (9.0.1)] or [HW]) we immediately get
the mod $p$ form of (1.23) from the above.

The proof of Theorem 1.3 is now complete. \qed

\Ack. The author would like to thank the referee for helpful comments.

\enddocument